\documentclass[10pt]{amsart}
\usepackage{amsmath}
\usepackage{enumerate}
\usepackage{pst-node}
\usepackage{tikz-cd}
\usepackage{amsopn}
\usepackage{amsmath, amsthm, amssymb, amscd, amsfonts, xcolor}
\usepackage[latin1]{inputenc}
\usepackage{graphicx}

\usepackage{euscript}

\usepackage{xcolor}
\usepackage{enumitem}
\usepackage{multirow}

\input xy
\xyoption{all}

\usepackage[OT2,OT1]{fontenc}
\newcommand\cyr{%
\renewcommand\rmdefault{wncyr}%
\renewcommand\sfdefault{wncyss}%
\renewcommand\encodingdefault{OT2}%
\normalfont
\selectfont}
\DeclareTextFontCommand{\textcyr}{\cyr}

\theoremstyle{plain}
\newtheorem{mainteo}{Theorem}

\newtheorem{teo}{Theorem}[section]

\theoremstyle{definition}

\numberwithin{equation}{section}

\def\Hol{\mathrm{H}}
\def\Ort{\mathrm{O}(n)}
\def\Nor[#1]#2{\mathrm{N}_{#1}(#2)}

\title{The holonomy group of a locally symmetric space}

\author{Antonio J. Di Scala}

\address{Dipartimento di Scienze Matematiche \it G.L. Lagrange \rm , Politecnico di Torino, Corso Duca degli Abruzzi 24, 10129, Torino, Italy.}

\email{antonio.discala@polito.it}
\thanks{Antonio J. Di Scala is a member of the research group in cryptography and
number theory (CrypTO) at Politecnico di Torino and GNSAGA at INdAM}


\begin{document}

\begin{abstract}
We show that the holonomy group of a connected Riemannian locally symmetric space (not necessarily complete) without local flat factor is compact and has finite index in its normalizer in the orthogonal group.
\end{abstract}

\maketitle

\section{Introduction}

The goal of this note is to provide a detailed proof of the following fact:

\begin{mainteo}\label{main} Let $(M^n,g)$ be a connected Riemannian locally symmetric space, not necessarily complete, without local de Rham flat factor.
Let $\Hol$ be its holonomy group and let $\Nor[\Ort]{\Hol}$ be its normalizer inside the orthogonal group $\Ort$. Then $\Hol$ and $\Nor[\Ort]{\Hol}$ are both compact, have the same dimension, and $\Hol$ has finite index inside $\Nor[\Ort]{\Hol}$.
\end{mainteo}

By {\it without local de Rham flat factor} we mean that the local de Rham decomposition, around each point of $M$ has no flat Riemannian factor or, equivalently,
the local holonomy group $\Hol^{Loc}_p$, acting on the tangent space $T_pM$, has no non-trivial fixed points (see Theorem \ref{locDe} below).\

Our motivation is a recent question posed by Graf and Patel in \cite{GP25}, related to an analogous claim in \cite[Prop. 10.114]{Bes87}.
Actually, our theorem is more general than the one of Besse. This proves that Besse's Prop. 10.114 holds even for general connected Riemannian locally symmetric spaces not necessarily complete.

\section{Preliminaries}

\subsection{Basic facts} To prove Theorem \ref{main} we need the following facts of Lie groups and Riemannian geometry:\

\vspace{.2cm}

\begin{itemize}

\item[Fact 1)] Let $\mathrm{G} \subset \Ort$ be a compact subgroup. Then its normalizer $\Nor[\Ort]{\mathrm{G}}$ is also compact.\

\vspace{.2cm}

\item[Fact 2)] A compact Lie group $\mathrm{G}$ has a finite number of connected components. Indeed, the quotient $\mathrm{G}/\mathrm{G}^0$ is a compact discrete Lie group hence finite.\

\vspace{.2cm}

\item[Fact 3)] A connected subgroup $\mathrm{G} \subset \Ort$ acting irreducibly on $\mathbb{R}^n$ is compact. This is \cite[Theorem 2, Appendix 5]{KN63}.\

\vspace{.2cm}

\item[Fact 4)] Let $(M,g)$ be a Riemannian manifold and let $\Hol^{Loc}_p$ be the local holonomy group at $p \in M$. Then $\Hol^{Loc}_p$ is a subgroup of the restricted holonomy group $\Hol_p^0$ at $p$.
If the metric $g$ is analytic then $\Hol^{Loc}_p = \Hol_p^0$. This is a direct consequence \cite[Proposition 10.1, page 95]{KN63} and \cite[Theorem 10.8, page 101]{KN63}.\

\vspace{.2cm}

\item[Fact 5)] Let $\mathrm{H}$ be a connected Lie subgroup of $\Ort$, acting on $\mathbb{R}^n$ as an s-representation and let $\Nor[\Ort]{\mathrm{H}}^0$ be the connected component of the normalizer of
$\mathrm{H}$ in $\Ort$. Then $\mathrm{H} = \Nor[\Ort]{\mathrm{H}}^0$. This is \cite[Lemma 6.2.2, page 192]{BCO03} or \cite[Lemma 5.2.2, page 184]{BCO16}.\

\end{itemize}

\vspace{.2cm}

Let me recall, following \cite[page 45]{BCO16} or \cite[page 46]{BCO03},  the definition of an s-representation.
Let $M = \mathrm{G}/\mathrm{K}$ be simply connected semisimple Riemannian symmetric space with $\mathrm{G} = \mathrm{Iso}^0(M)$ semisimple and
$\mathrm{K} = \mathrm{G}_p$ the isotropy group at $p \in M$. Since $M$ is simply connected and $\mathrm{G}$ is connected, $\mathrm{K}$ is also connected.
The isotropy representation of $\mathrm{G}/\mathrm{K}$
at $p$, that is, the Lie group homomorphism \[ \rho: K \to \mathrm{SO}(T_p M) \, \, \, \, , \, \, \, \, k \to (\mathrm{d}k)_p , \]
is called an s-representation.
We say that a connected Lie subgroup $\mathrm{H}$ of $\Ort$ is acting on $\mathbb{R}^n$ as an s-representation if there is a simply connected semisimple Riemannian symmetric space $\mathrm{G}/\mathrm{K}$
and isomorphisms $\xi : \mathrm{H} \to \mathrm{K}$ and $\iota: \mathrm{SO}(n) \to \mathrm{SO}(T_pM) $ making commutative the following diagram:

\[
\begin{tikzcd}
\mathrm{H} \arrow[r, "\xi"] \arrow[d, "i"']
& \mathrm{K} \arrow[d, "\rho"] \\
\mathrm{SO}(n) \arrow[r, "\iota"']
& \mathrm{SO}(T_pM)
\end{tikzcd}
\]

where $i: \mathrm{H} \to \mathrm{SO}(n)$ is the homomorphism given by the inclusion of the connected subgroup $\mathrm{H}$ of $\Ort$.

\subsection{Local de Rham decomposition}  Here we recall the local de Rham decomposition around a point $p \in M$ of a Riemannian manifold $(M,g)$ and the corresponding factorization of the local holonomy group $\Hol^{Loc}_p$.
We will call Local de Rham decomposition the following theorem:

\vspace{.2cm}

\begin{teo} \label{locDe} Let $(M,g)$ be a Riemannian manifold and let $\Hol^{Loc}_p \subset \mathrm{SO}(T_pM)$ the local holonomy group at $p$. Then the tangent space $T_pM$ splits as orthogonal direct sum of $\Hol^{Loc}_p$-invariant subspaces:
\[ T_pM = V_0 \oplus V_1 \oplus \cdots \oplus V_k \]
where $V_0$ is the maximal linear subspace of $T_pM$ on which $\Hol^{Loc}_p$ acts trivially and on each subspace $V_j$,$(0 < j \leq k)$, the $\Hol^{Loc}_p$-action is irreducible.

Moreover, by parallel transport each $V_j$ induces a parallel distribution around $p$ whose integral manifolds are denoted as $\mathcal{V}_j$.
Let $g_j$ be the restriction of the Riemannian metric $g$ to $\mathcal{V}_j$.  Then there is a neighborhood $\mathcal{V}$ of $p$ isometric to the product $\prod_{j=0}^k (\mathcal{V}_j , g_j)$.
The local holonomy group $\Hol^{Loc}_p$ also splits as a product of connected subgroups $\Phi_j$, $(0 < j \leq k)$ :
\[ \Hol^{Loc}_p = \prod_{j=1}^k \Phi_j\]
where $\Phi_j$ acts irreducibly on $V_j$ and trivially on $V_i$ for $i \neq j$.\\

Each $\Phi_j$ is compact and hence $\Hol^{Loc}_p$ and $\Nor[\Ort]{\Hol^{Loc}_p}$ are compact Lie groups.
\end{teo}

The above is proved in \cite[Theorem 5.4., page 185]{KN63}. The Riemannian manifold $(\mathcal{U}_0,g_0)$ is flat and called {\it the flat local de Rham factor}.
Observe that by Fact 3) each $\Phi_j$ is compact hence $\Hol^{Loc}_p$ is compact. Then by Fact 1) $\Nor[\Ort]{\Hol^{Loc}_p}$ is also compact.

\section{Proof of Theorem \ref{main}}

Let \((M,g)\) be a connected locally symmetric space without local de Rham flat factor, and fix \(p \in M\).
Let \(\Hol^{loc}_p\) be its local holonomy group at \(p\).
By the local de Rham theorem (Theorem \ref{locDe}), the tangent space decomposes as
\[
T_pM = V_1 \oplus \cdots \oplus V_k ,
\]
where each \(V_j\) is an irreducible \(\Hol^{loc}_p\)-module and the group \(\Hol^{loc}_p\) splits as a product
\[
\Hol^{loc}_p = \Phi_1 \times \cdots \times \Phi_k ,
\]
with \(\Phi_j\) acting irreducibly on \(V_j\) and trivially on the other factors.

Because \((M,g)\) is locally symmetric, each \((\mathcal{V}_j,g_j)\) in the local decomposition is locally isometric to an open subset of an irreducible simply connected symmetric space \(M_j\).
The metric of a locally symmetric space is real analytic, hence by Fact 4 we have \(\Hol^{loc}_p = \Hol_p^0\) (the restricted holonomy group).
Consequently each \(\Phi_j\) is (isomorphic to) the holonomy group of the irreducible symmetric space \(M_j\); therefore \(\Phi_j\), and hence \(\Hol^{loc}_p\), acts as an s-representation.

Applying Fact 5 to \(\mathrm{H} = \Hol^{loc}_p\) gives
\[
\Hol^{loc}_p = \bigl(\mathrm{N}_{\mathrm{O}(T_pM)}(\Hol^{loc}_p)\bigr)^0 .
\]

Now let \(\Hol\) be the full holonomy group of \((M,g)\) and \(\Hol^0\) its identity component.
From Fact 4 we have \(\Hol^0 = \Hol^{loc}_p\).
Thus \(\Hol^0 = \bigl(\mathrm{N}_{\mathrm{O}(T_pM)}(\Hol^0)\bigr)^0\).
In particular \(\Hol^0\) is compact (Fact 3), and so is \(\mathrm{N}_{\mathrm{O}(T_pM)}(\Hol^0)\) (Fact 1).
Since \(\Hol \subset \mathrm{N}_{\mathrm{O}(T_pM)}(\Hol^0)\), the group \(\Hol\) is also compact; another application of Fact 1 shows that \(\mathrm{N}_{\mathrm{O}(T_pM)}(\Hol)\) is compact as well.

Finally, because \(\Hol^0 = \bigl(\mathrm{N}_{\mathrm{O}(T_pM)}(\Hol)\bigr)^0\) (the inclusion \(\Hol^0 \subset \bigl(\mathrm{N}_{\mathrm{O}(T_pM)}(\Hol)\bigr)^0 \subset \mathrm{N}_{\mathrm{O}(T_pM)}(\Hol^0)\) forces equality), the two groups \(\Hol\) and \(\mathrm{N}_{\mathrm{O}(T_pM)}(\Hol)\) have the same dimension.
Moreover, \(\mathrm{N}_{\mathrm{O}(T_pM)}(\Hol)\) has finitely many connected components (Fact 2), so the index of \(\Hol\) in its normalizer is finite.
This completes the proof.



\begin{thebibliography}{BCO00}

\bibitem[BCO03]{BCO03}
Berndt, J. ; Console, S. and Olmos, C.:
\textit{Submanifolds and holonomy}
Chapman \& Hall/CRC , Research Notes in Mathematics
434 (2003).


\bibitem[BCO16]{BCO16}
Berndt, J. ; Console, S. and Olmos, C.:
\textit{Submanifolds and holonomy} Second Edition,
Chapman \& Hall/CRC Research Notes in Mathematics, 434,  2016.

\bibitem[Bes87]{Bes87} Besse, A.L.:
\textit {Einstein Manifolds}
Vol. 10 in the series: Ergebnisse der Mathematik und ihrer Grenzgebiete.
Springer-Verlag Berlin 1987.


\bibitem[GP25]{GP25} Graf, P. and Patel, A.:
\textit {Slope zero tensors, uniformizing variations of Hodge structure and quotients of tube domains}
https://arxiv.org/pdf/2510.15039


\bibitem[KN63]{KN63} Kobayashi, S. and Nomizu, K.:
\textit {Foundations of differential geometry}
Vol I, Interscience Publishers, (1963).



\end{thebibliography}
\end{document}